\input amstex 
\documentstyle{amsppt}
\input bull-ppt
\keyedby{bull302/lbd}

\topmatter
\cvol{27}
\cvolyear{1992}
\cmonth{October}
\cyear{1992}
\cvolno{2}
\cpgs{243-245}
\title The 1-, 2-, and 3-characters determine a 
group\endtitle
\author H.-J. Hoehnke and K. W. Johnson\endauthor
\address Mendelstr. 4, 0-1100 Berlin, Germany and 
Brandenburgische
Landesuniversit\" at, Fachbereich Mathematik, AM Neuen 
Palais, 0-1571
Potsdam-Sanssouci, Germany\endaddress
\address Department of Mathematics,
The Pennsylvania State University Abington,
Pennsylvania 19001\endaddress
\ml kwj1 \@ psuvm.bitnet\endml
\date March 20, 1991 and, in revised form, February 20, 
1992\enddate
\abstract A set of invariants for a finite group is 
described. These arise
naturally from Frobenius' early work on the group 
determinant and provide an
answer to a question of Brauer. Whereas it is well known 
that the ordinary
character table of a group does not determine the group 
uniquely, it is a
consequence of the results presented here that a group is 
determined uniquely
by its ``3-character'' table.\endabstract
\endtopmatter

\document
\heading 1. Introduction\endheading
\par
Given a finite group $G$ of order $n$ its {\it group 
matrix\/} $X_G$ is defined
as follows. Let $x_e,x_{g_2},\dotsc,x_{g_n}$ be variables 
indexed by the
elements of $G$. The $n\times n$ matrix $X_G$ is defined 
to be the matrix whose
$(i,j)\roman{th}$ entry is $x_{g_ig_j^{-1}}$. The {\it 
group determinant\/}
$\Theta_G$ of $G$ is $\det(X_G)$.
\par
Although the group determinant first appeared in 
Frobenius' 1896 paper
\cite{Fr1} in which he introduced characters for arbitrary 
groups, the roots of
the work go back to the discussions of Gauss on the 
composition of equivalence
classes of quadratic forms. There is a good discussion of 
the background to
Frobenius' work in [Ha1, Ha2]. A natural question arises 
as to whether
nonisomorphic groups necessarily have distinct group 
determinants. This was
posed by one of the authors in 1986 and was answered by 
Formanek and Sibley
\cite{FS} in 1990, the positive answer being somewhat 
surprising. In fact there
were already two early remarks in the literature apart 
from the work of Gauss
quoted above that pointed in this direction. Both of these 
were described in
apostrophies as strange (``merkw\"urdig'') by their authors 
and were disregarded
for a long time within the development of algebra. The 
first was by Frobenius
in \cite{Fr2} on matrix transformations and the second 
occurred in \cite{B} on
the subject of maximal orders of quaternion algebras. 
Today these phenomena are
subsumed under the thesis ``from norm invariance to 
constructive structure
theory and (noncommutative) arithmetics.''
\par
A further look at Frobenius' early work reveals functions 
$\chi^{(k)}\colon
G^k\to\Bbb C$ defined below in the case where $k\le3$ that 
appeared in his
algorithm to calculate the factor of a group determinant 
that corresponds to an
irreducible character $\chi$. \pagebreak
Some of the properties of 
these
``$k$-characters'' have been explored in \cite{J} where it 
is shown that the
$k$-characters corresponding to distinct irreducible 
characters $\chi_i,\chi_j$
are orthogonal in the sense that
$$\sum_{\underline g\in G^k}\chi_i^{(k)}(\underline
g)\overline{\chi}_j^{(k)}(\underline g)=0.$$
There is also given in \cite{J} a 2-character table of a 
group $G$ that
consists of the 2-characters corresponding to the 
irreducible representations
of degrees greater than 1 together with ``degenerate'' 
characters.
\par
One of the questions raised by Brauer in \cite{Br} is that 
of determining what
extra information may be added to the (ordinary) character 
table to completely
determine a group. As a consequence of Frobenius' work and 
the Formanek-Sibley
result, the knowledge of the $k$-characters $\chi^{(k)}$ 
of a group for all $k$
and all irreducible characters $\chi$ is sufficient to 
determine a group. Note
that if $k>\deg(\chi)$ then $\chi^{(k)}=0$. The 
$k$-characters thus provide an
answer to Brauer's question, but since the amount of work 
involved in
calculating the $k$-characters for large $k$ is 
prohibitively large, it becomes
interesting to examine the question of the extent to which 
the knowledge of the
$k$-characters for small values of $k$ determine a group.
\heading 2. 3-characters\endheading
\par
The authors would like to announce the following result.
\thm{Theorem} Let $G$ be a finite group, with irreducible 
characters
$\chi_1,\dotsc,\chi_r$ over an algebraically closed field 
$K$ with
$\roman{char}\,K\not=2$ and $\roman{char}\,K\nmid|G|$. 
Then $G$ is determined
up to isomorphism by the $\chi_j^{(k)}$, $j=1,\dotsc,r$, 
$k=1,2,3$.
\ethm
\par
Explicit definitions of the 1-, 2-, and 3-characters 
corresponding to a
character $\chi$ are as follows:
$$\gather
\chi^1(g)=\chi(g),\qquad 
\chi^2(g,h)=\chi(g)\chi(h)-\chi(gh),\\
{\aligned\chi^3(g,h,m)=&\chi(g)\chi(h)\chi(m)-\chi(g)%
\chi(hm)-\chi(h)\chi(gm)\\
&-\chi(m)\chi(gh)+\chi(ghm)+
\chi(gmh).\endaligned}\endgather$$
\par
The proof of the theorem depends on results in [Ho1, Ho2]; 
the following is a
brief outline:
\par
(i) To each irreducible character $\chi$ of a finite group 
$G$ there is
associated a factor $\varphi_\chi$ of $\Theta_G$.
\par
(ii) Any factor $\varphi$ of $\Theta_G$ is a norm-type 
form. In particular,
$\varphi(xy)=\varphi(x)\varphi(y)$ where $x$ and $y$ are 
generic elements of
the group algebra of $G$.
\par
(iii) Let $A$ be a finite-dimensional algebra over the 
field $K$, and let
$\{\omega_1,\dotsc,
\omega_n\}$ be a basis. Define the structure constants
$\{w_{ij}^k\}$ by
$$\omega_i\omega_j=\sum w_{ij}^k\omega_k.$$
Let $N$ be a norm-type form on $A$, with
$$N(\lambda-x)=\lambda^m-s_1(x)\lambda^{m-1}+\dotsb+
(-1)^ms_m(x)$$
where $\lambda$ is an indeterminate and $x$ is a generic 
element,
$$x=x_1\omega_1+\dotsb+x_n\omega_n.$$
If the discriminant of $N$ is nonzero (which is true for 
example when $N$ is
the group determinant) it follows that from the knowledge 
of $s_1(x)$,
$s_2(x)$, and $s_3(x)$ the ``symmetrised'' structure 
constants
$$w_{(ij)}^k=w_{ij}^k+w_{ji}^k$$
can be determined.
\par
(iv) By [Bo, \S4, Exercise 26] if $f\colon G\to H$ is a 
bijection between finite
groups such that $f(gh)$ is either $f(g)f(h)$ or 
$f(h)f(g)$, then $f$ is either
an isomorphism or an anti-isomorphism, and hence $G$ and 
$H$ are isomorphic.
\par
(v) It follows from (iii) and (iv) that the 3-character of 
the regular
representation (which is essentially the same as $s_3(x)$ 
above when
$N=\Theta_G)$ determines $G$.
\par
(vi) The 3-character of the regular representation of $G$ 
can be formed from
the 1-, 2-, and 3-characters associated to the irreducible 
representations of
$G$.
\par
This provides, in some sense, a more satisfactory set of 
invariants for a
finite group. It has also been possible to give a 
constructive proof of the
Formanek-Sibley theorem using the ideas in the proof of 
the above theorem.
\par
Another set of invariants for finite groups has been given 
by Roitman in
\cite{Ro}. These invariants consist of an apparently 
infinite set of integers
defined for pairs of integers $n,k$ that are calculated as 
the coefficient of
the identity in ô¨the $k\roman{th}$ powers of certain 
elements of $\Bbb
ZG\otimes\dotsb\otimes\Bbb ZG$ $(n$ factors), identifying 
this ring with $\Bbb
ZG^n$. We refer the reader to \cite{Ro} for the details, 
which are of a
technical nature. The 1-, 2-, and 3-characters associated 
to the irreducible
representations of the group seem to be much more 
accessible, and moreover, a
3-character table can be constructed that has convenient 
orthogonality
properties. The question of whether the 1- and 
2-characters alone determine a
group is addressed in \cite{JS} where it is shown that 
there exist
nonisomorphic pairs of groups with the same 2-character 
table, i.e., with the
same 1- and 2-characters.


\Refs
\ra\key{Ho2}
\ref
\key B
\by H. Brandt
\paper Idealtheorie in Quaternionenalgebren
\jour Math. Ann.
\vol 99
\yr 1928
\pages 1--29
\endref
\ref
\key Bo
\by N. Bourbaki
\book Alg\`ebre. {\rm I}
\publaddr Paris
\yr 1970
\endref
\ref
\key Br
\by R. Brauer
\book Representations of finite groups
\bookinfo Lectures in Modern Mathematics {\rm(T. L. Saaty, 
ed.)}
\publaddr Wileÿ\yr 19øÄ63
\pages 133--175
\endref
\ref
\key FS
\by E. Formanek and D. Sibley
\paper The group determinant determines the group
\jour Proc. Amer. Math. Soc. {\bf 112} (1991), 649--656
\endref
\ref
\key Fr1
\by G. Frobenius
\paper \"Uber die Primfaktoren der Gruppendeterminante
\jour S'ber. Akad. Wiss. Berlin
\yr 1896
\pages 1343--1382
\endref
\ref
\key Fr2
\bysame
\paper \"Uber die Darstellung der endlichen Gruppen durch 
lineare
Substitutionen
\jour S'ber. Akad. Wiss. Berlin
\yr 1898
\pages 944--1015
\endref
\ref
\key Ha1
\by T. Hawkins
\paper The origins of the theory of group characters
\jour Arch. Hist. Exact Sci.
\vol 7
\yr 1971
\pages 142--170
\endref
\ref
\key Ha2
\bysame
\paper New light on Frobenius' creation of the theory of 
group characters
\jour Arch. Hist. Exact Sci.
\vol 12
\yr 1974
\pages 217--243
\endref
\ref
\key Ho1
\by H.-J. Hoehnke
\paper \"Uber komponierbare Formen und konkordante 
hyperkomplexe Gr\"ossen
\jour Math. Z.
\vol 70
\yr1958
\pages 1--12
\endref
\ref
\key Ho2
\bysame
\paper \"Uber Beziehungen zwischen Problemen von H. Brandt 
aus der Theorie
der Algebren und den Automorphismen der Normenform
\jour Math. Nachr.
\vol 34
\yr 1967
\pages 229--255
\endref
\ref
\key J
\by K. W. Johnson
\paper On the group determinant
\jour Math. Proc. Cambridge Philos. Soc.
\vol 109
\yr 1991
\pages 299--311
\endref
\ref
\key JS
\by K. W. Johnson and Surinder Sehgal
\book The {\rm2}-character table does not determine a 
group, {\rm preprint}
\endref
\ref
\key Ro
\by M. Roitman
\paper A complete set of invariants for finite groups
\jour Adv. in Math.
\vol 41
\yr 1981
\pages 301--311
\endref

\endRefs
\enddocument